\documentclass[
final
, nomarks
]{dmtcs-episciences}

\usepackage[utf8]{inputenc}
\usepackage{subfigure}

\usepackage{diagbox}
\usepackage[table]{xcolor}
\usepackage{amsmath}

\def\N{\mathbb{N}}
\def\Av{\textsc{Av}}
\def\pop{\textsc{pop}}

\newtheorem{theorem}{Theorem}[section]
\newtheorem{lemma}[theorem]{Lemma}
\newtheorem{corollary}[theorem]{Corollary}
\newtheorem{example}[theorem]{Example}
\newtheorem{proposition}[theorem]{Proposition}
\newtheorem{conjecture}[theorem]{Conjecture}

\usepackage[round]{natbib}

\author[Nathana\"el Hassler and Sergey Kirgizov]{Nathana\"el Hassler
  \and Sergey Kirgizov}
\title[Emerging consecutive pattern avoidance]{Emerging consecutive pattern avoidance}
\affiliation{
  Université Bourgogne Europe, LIB UR 7534, 21000 Dijon, France
}
\begin{document}
% This is only used if you are compiling for a volume before vol 25
% \publicationdetails{VOL}{2015}{ISS}{NUM}{SUBM}
% This is the new form of collecting the data, starting with vol 25
\publicationdata{vol. 28:1, Permutation Patterns
2025}{2026}{9}{10.46298/dmtcs.16861}{2025-11-05; 2025-11-05;
2026-10-14}{2026-10-30}
%% \publicationdata
%% {vol. 25:3 special issue for main purpose}
%% {2022}
%% {1}
%% {10.46298/dmtcs.10472}
%% %{1998-10-14; 1998-10-14; 2002-07-19; 2014-02-05; 2015-09-09; 2022-12-25}
%% %{2022-12-3}
%% {2022-12-3; None}
%% {2023-1-1}

\maketitle
\begin{abstract}
  In this note we study the \emph{asymptotic popularity}, that is, the
  limit probability to find a given consecutive pattern at a random
  position in a random permutation in the eighteen classes of
  permutations avoiding at least two length 3 consecutive patterns. We
  show that for ten classes, this popularity can be readily deduced
  from the structure of the permutations.  By combining analytical and
  bijective approaches, we study in detail two more involved
  cases. The problem remains open for five classes.
\end{abstract}

\section{Introduction and notation}

We write a permutation $\pi\in\mathcal{S}_n$ as a word $\pi=a_1\ldots
a_n$ whose letters are $\{a_1,\ldots,a_n\}=\{1,\ldots,n\}$. A
\textit{pattern} $p$ of length $r$ is an element of
$\mathcal{S}_r$. We usually say that $\pi$ contains an occurrence of
the pattern $p$ if there exists a subsequence $1\leq i_1<\ldots
<i_r\leq n$ such that $a_{i_1}\ldots a_{i_r}$ is order-isomorphic to
$p$. In this note we focus on \textit{consecutive} patterns. We say
that $\pi$ contains a consecutive occurrence of the pattern
$p$ if there exists a subsequence of consecutive letters
$a_{i}a_{i+1}\ldots a_{i+r-1}$ of $\pi$ that is order-isomorphic to
$p$. We say that $\pi$ \textit{avoids} the consecutive pattern $p$ if
it does not contain any consecutive occurrence of $p$.

\medskip

\cite{kit}, and then~\cite{kitmans,kitmans2}, presented the
enumeration of classes of permutations avoiding at least two length 3
consecutive patterns.  In this work, we study the \emph{asymptotic
popularity}, that is, the limit probability to find a given pattern of
length 3 at a random position in a random permutations in the eighteen
avoidance classes from Kitaev's and Mansour's works. We show that, in
certain cases, some of the remaining patterns disappear
asymptotically. It is a quite enchanting fact.

\medskip

In the wide realm of interesting papers on permutation patterns, we
would like to highlight those most relevant to our
work. \cite{bona1,bona2} and \cite{homberger} studied the same problem
for classical patterns (i.e. non-consecutive) in classes of
permutations avoiding one length 3 pattern. In particular, they proved
that among the permutations avoiding the classical pattern $123$, and
the ones avoiding $132$, the decreasing pattern $321$ appears
asymptotically with probability 1, while the four other length 3
patterns disappear.

\cite{jan17, jan19} considered limit laws for the distributions of
classical patterns of length 3 in permutations avoiding classical
patterns 132 and 321. \cite{borga1} studied the local limit of
permutations avoiding one classical pattern of length 3, which is the
same as studying the proportion of every consecutive patterns, of any
length. In particular, he proved that in the class of permutations
avoiding 321, consecutive patterns whose inverse has more than one
descent disappear asymptotically. In another paper,~\cite{borga}
introduced a method based on generating trees to study asymptotic
normality of consecutive pattern occurrences in permutations avoiding
certain classical patterns.

\cite{elizalde-noy} presented a method based on increasing trees and
the boxed product (see~\cite{flaj}) to study distributions and
avoidance of certain consecutive patterns in permutations.  The same
authors, in another paper (see~\cite{elizalde-noy-clusters}), showed
how the Goulden-Jackson cluster method (see~\cite{gj1, gj2}) can be
adapted to enumerate permutations that avoid consecutive patterns.

\cite{barnabei} studied joint distributions of consecutive patterns of
length 3 in the set of permutations avoiding the classical pattern
312.  They did this by observing how the patterns are transformed by
Krattenthaler's bijection (see~\cite{krattenthaler}) between such
permutations and Dyck paths, and how Deutsch's involution
(see~\cite{deutsch}) on Dyck paths helps this process. Their pattern
transfer method is similar to what Baril, Burstein and Kirgizov did in
their work about faro words and permutations (see~\cite{faro}).  Their
paper is a precursor to the article you are holding in your hands.

\medskip

Now let us introduce some bijections from $\mathcal{S}_n$ to $\mathcal{S}_n$, that preserve the occurrences of consecutive patterns. They will allow us to reduce significantly the number of multi-avoidance classes to study. The \textit{reverse} $R(\pi)$ of a permutation $\pi=a_1\ldots a_n$ is
the permutation $a_n\ldots a_1$. The \textit{complement} $C(\pi)$ is
the permutation $(n+1-a_1)\ldots (n+1-a_n)$. We denote by $R\circ C$ the
composition of $R$ and $C$. For example, $R(35214)=41253$,
$C(35214)=31452$ and $(R\circ C)(35214)=25413$. These bijections
preserve the occurrences of consecutive patterns; indeed for
$T\in\{R,C,R\circ C\}$, $\pi\in\mathcal{S}_n$, and a consecutive
pattern $p$, $\pi$ has an occurrence of $p$ if and only if $T(\pi)$
has an occurrence of $T(p)$. As explained by~\cite{kit}, this enables
us to focus on equivalence classes of the avoidance classes under
the action of these 3 bijections $R,C$ and $R\circ C$.

Table~\ref{tab:results} summarizes our results, presenting the
eighteen classes as they appear in the work by~\cite{kit}. In this
table, an empty cell means that corresponding patterns are avoided
\emph{by design}, while 0 denotes that the respective pattern
disappears asymptotically (the probability to find this pattern at a
random position in a random permutation tends to 0, as the permutation
length grows).  The values 1/2 and 1/4 present in this table should
also be understood in the asymptotic sense.  There are two ``N/A'' for
Class 4 because this class is empty for $n>3$.  Question marks
indicate open problems.

\begin{table}[htb]
  \begin{center}
    \begin{tabular}{|l|c|c|c|c|c|c|}
    \hline
        \diagbox{Class}{Pattern} & $123$ & $132$ & $213$ & $231$ & $312$ & $321$\\
         \hline
        $\mathbf{1}$ (simple, Sec.~\ref{sec:simple}) & \cellcolor{gray!20} & \cellcolor{gray!20} & 1/2 & 1/2 & \cellcolor{gray!20} & \cellcolor{gray!20} \\
        \hline
        $\mathbf{2}$ (simple, Sec.~\ref{sec:simple})& \cellcolor{gray!20} & \cellcolor{gray!20} & \cellcolor{gray!20} & 0 & \cellcolor{gray!20} & 1\\
        \hline
        $\mathbf{3}$ (simple, Sec.~\ref{sec:simple}) & 1/2 & \cellcolor{gray!20} & \cellcolor{gray!20} & \cellcolor{gray!20} & \cellcolor{gray!20} & 1/2\\
         \hline
        $\mathbf{4}$ (simple, Sec.~\ref{sec:simple}) & \cellcolor{gray!20} & \cellcolor{gray!20} & N/A & \cellcolor{gray!20} & N/A & \cellcolor{gray!20}\\
         \hline
        $\mathbf{5}$ (simple, Sec.~\ref{sec:simple}) & 1 & \cellcolor{gray!20} & \cellcolor{gray!20} & 0 & \cellcolor{gray!20} & \cellcolor{gray!20}\\
         \hline
        $\mathbf{6}$ (simple, Sec.~\ref{sec:simple}) & \cellcolor{gray!20} & \cellcolor{gray!20} & \cellcolor{gray!20} & 1/2 & 1/2 & \cellcolor{gray!20}\\
         \hline
        $\mathbf{7}$ (done by~\cite{faro}) & \cellcolor{gray!20} & \cellcolor{gray!20} & \cellcolor{gray!20} & 1/2 & 1/2 & 0\\
         \hline
         $\mathbf{8}$ (simple, Sec.~\ref{sec:simple}) & \cellcolor{gray!20} & \cellcolor{gray!20} & 0 & \cellcolor{gray!20} & 0 & 1\\
         \hline
        $\mathbf{9}$ (simple, Sec.~\ref{sec:simple}) & 1/2 & \cellcolor{gray!20} & \cellcolor{gray!20} & \cellcolor{gray!20} & 0 & 1/2\\
         \hline
        $\mathbf{10}$ (open)& \cellcolor{gray!20} & \cellcolor{gray!20} & ? & ? & \cellcolor{gray!20} & ?\\
         \hline
        $\mathbf{11}$ (Section~\ref{sec:class11}) & \cellcolor{gray!20} & \cellcolor{gray!20} & 1/4 & 1/2 & 1/4 & \cellcolor{gray!20}\\
         \hline
        $\mathbf{12}$ (open)& \cellcolor{gray!20} & ? & ? & \cellcolor{gray!20} & \cellcolor{gray!20} & ?\\
         \hline
        $\mathbf{13}$ (open)& \cellcolor{gray!20} & ? & ? & \cellcolor{gray!20} & ? & ?\\
         \hline
        $\mathbf{14}$ (open)& ? & ? & \cellcolor{gray!20} & \cellcolor{gray!20} & ? & ?\\
         \hline
        $\mathbf{15}$ (open) & ? & \cellcolor{gray!20} & \cellcolor{gray!20} & ? & ? & ?\\
         \hline
        $\mathbf{16}$ (simple, Sec.~\ref{sec:simple})& \cellcolor{gray!20} & 1/4 & 1/4 & 1/4 & 1/4 & \cellcolor{gray!20} \\
         \hline
        $\mathbf{17}$ (Section~\ref{sec:class17})& \cellcolor{gray!20} & \cellcolor{gray!20} & 1/4 & 1/2 & 1/4 & 0\\
         \hline
        $\mathbf{18}$ (simple, Sec.~\ref{sec:simple}) & 1/2 & \cellcolor{gray!20} & 0 & \cellcolor{gray!20} & 0 & 1/2\\
         \hline
    \end{tabular}\vspace{1ex}
    
    \caption{The asymptotic popularity patterns among eighteen avoidance classes.}
    \label{tab:results}
  \end{center}
\end{table}

For consecutive patterns $p_1,\ldots,p_k$ and $n\in\N$, we denote by
$\Av_n(p_1,\ldots,p_k)$ the set of permutations of length $n$ that avoid
all of the consecutive patterns $p_1,\ldots,p_k$, and by
$\Av(p_1,\ldots,p_k)$ the set of all permutations that avoid
$p_1,\ldots,p_k$.

Let $\mathcal{A}_n:=\Av_n(p_1,\ldots,p_k)$. For a pattern
$p\not\in\{p_1,\ldots,p_k\}$, we denote by
$\mathbf{p}_n^{\mathcal{A}}$ the \textit{popularity} of the pattern
$p$ in the class $\mathcal{A}_n$, that is, the total number of
occurrences of $p$ in $\mathcal{A}_n$. Now we define the
\textit{asymptotic popularity} of $p$ in the class $\mathcal{A}:=\Av(p_1,\ldots,p_k)$ by

\begin{equation}
\pop_{\mathcal{A}}(p):=\lim_{n\to\infty}\frac{\mathbf{p}_n^{\mathcal{A}}}{n|\mathcal{A}_n|},
\label{eq:limpop}
\end{equation}
when the limit exists.  We will use $\pop_k$, where $k$ refers to
Class $k$ from Table~\ref{tab:results}.  When the class is clear from the
context, we simply write $\pop$.

The note is organized as follows: in Section~\ref{sec:simple} we
present some classes for which asymptotic popularities are easily
obtained. In Sections~\ref{sec:class11} and~\ref{sec:class17} we
compute these popularities for two more classes which are more complex. The methods
combine analytic arguments regarding the generating functions and a
bijective argument connecting the permutations from the class to
involutions. Finally, in Section~\ref{sec:open} we offer a conjecture
and formulate several questions for future research.

\section{Simple classes}
\label{sec:simple}

 In this section, we analyze the classes for which the asymptotic
 popularities are easily computable: classes
 1,\,2,\,3,\,5,\,6,\,8,\,9,\,16 and 18 (see Table~\ref{tab:results}). Note that the
 problem is not interesting for Class 4, $\Av_n (123, 132, 231, 321)$, since the class is empty
 as soon as $n>3$. Class~7,
 $\Av_n(123,132,213)$, has been done by
 Baril, Burstein and Kirgizov (see Remark 3.8,~\cite{faro}), using a
 bijection between such permutations and dispersed Dyck paths that
 transfers patterns between these two sets of combinatorial objects in
 a nice way. Their result was the motivation
 for the present work.

\begin{figure}[ht]
  \begin{center}
  \includegraphics[width=0.6\textwidth]{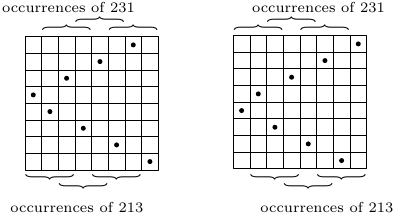}
  \caption{The only 2 permutations of Class 1, $\Av_n(123,132,312,321)$ for $n=8$.}
  \label{fig:class1}
  \end{center}
\end{figure}

Consider Class 1: $\Av_n(123,132,312,321),$ which contains only 2
permutations when $n>1$.  These permutations have a simple alternating
structure shown in Figure~\ref{fig:class1}. It it clear that the
asymptotic popularities of the two remaining patterns are equal: $\pop
(213) = \pop(231) = 1/2$.

Class 2, $\Av_n(123,132,213,312)$, also contains just 2 permutations
when $n>1$. The pattern $231$ may appear only once, at the very
beginning of a permutation. After that, we observe exclusively
occurrences of $321$ (for $n>3$), so $\pop (231) = 0$ and $\pop (321)
= 1$.

   Class 3, $\Av_n(132,213,231,312)$, for $n>1$ consists of two
   permutations: $123 \ldots n$ and $n (n-1) \ldots 321$, thus $\pop
   (123) = \pop (321) = 1/2$.

   Class 4, $\Av_n(123,132,231,321)$, is empty for $n>3$.

   A permutation from Class 5, $\Av_n(132,213,312,321)$, may have at
   most one occurrence of pattern $231$, which, when it exists, has to be in the last three elements of the permutation. So $\pop (231) = 0$ and
   $\pop(123) = 1$. There are $n-1$ permutations in Class 5, for
   $n>2$.

   For $n>3$, any permutation from Class 6, $\Av_n(123,132,213,321)$,
   is a sequence of overlaps of two alternating patterns $231$ and
   $312$.  We have $\pop (231) = \pop (312) = 1/2$.

   For $n>3$, Class 8, $\Av_n(123,132,231)$, consists of $n$
   permutations. Any permutation from this class starts with a
   sequence of descents. At the end it may have one occurrence of
   pattern $213$ or $312$. It follows that $\pop (213) = \pop (312) =
   0$ and $\pop (321) = 1$.

   A typical permutation from Class 9, $\Av_n(132,
   213, 231)$, starts with a sequence of descents,
   has one occurrence of pattern $312$, and ends with a sequence
   of ascents. The permutations $123 \ldots n$ and $n (n-1) \ldots 321$
   are also included. So, $\pop (123) = \pop (321)
   = 1/2$ and $\pop (312) = 0$.
 
Class 16, $\Av_n(123,321)$, can be solved directly
with a symmetry argument. Indeed, this class is stable under the
action of the 3 bijections $R,C$ and $R\circ C$. So any occurrence of
the pattern $132$ in this class is uniquely mapped to an
occurrence of $231$ in the same class through $R$, an
occurrence of $312$ through $C$, and an occurrence of
$213$ through $R\circ C$. Hence
$\pop_{16}(132)=\pop_{16}(231)=\pop_{16}(312)=\pop_{16}(213)=1/4$.

\medskip

Let us describe in detail the case of Class 18,
$\Av_n(132,231)$.

\begin{figure}[ht]
  \begin{center}
    \includegraphics[width=0.8\textwidth]{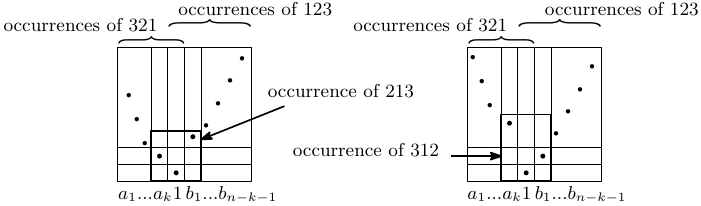}
    \caption{General structure of a permutation from $\Av_n(132,231)$.}
    \label{fig:class18}
  \end{center}
\end{figure}

As shown by~\cite{kit}, we have the permutations of the form
$a_1\ldots a_{k} 1 b_{1}\ldots b_{n-k-1}$, where $a_1\ldots a_{k}$ is
a decreasing sequence, and $b_{1}\ldots b_{n-k-1}$ an increasing
sequence (see Figure~\ref{fig:class18}). In such a permutation, there
are exactly $k-1$ occurrences of the pattern $321$ (except for $k=0$,
where there is no occurrence), so
\begin{equation}
  \mathbf{321}_n=\sum_{k=1}^{n-1}\binom{n-1}{k}(k-1)=(n-1)\cdot 2^{n-2}-2^{n-1}+1.
\end{equation}

Since length $n$ permutations in Class 18 are enumerated by $2^{n-1}$,
we have
\begin{equation}
\pop_{18}(321)=\underset{n\to\infty}{\lim}\frac{\mathbf{321}_n}{n\cdot
  2^{n-1}}=\frac{1}{2}.
\end{equation}

Since this class is reverse-invariant, we deduce that
$\pop_{18}(123)=\pop_{18}(321)=1/2$ and
$\pop_{18}(213)=\pop_{18}(312)=0$.

Another way to prove the same result is to notice that every
permutation in $\Av_n(132, 231)$ except identity and anti-identity
contains a single occurrence of $213$ or $312$ but not
both. Therefore, $\mathbf{213}_n + \mathbf{312}_n = |\Av_n(132, 231)|
- 2$, and we have $\pop_{18}(213)=\pop_{18}(312)=0$.
  
% Similarly we conclude that $\pop_{18}(123)=1/2$, and $\pop_{18}(213)=\pop_{18}(312)=0$.

\section{Avoiding $123$, $132$ and $321$}
\label{sec:class11}

 We consider the asymptotic relative popularity of the patterns 213,
 231 and 312 in $A_n:=\Av_n(123,132,321)$ (Class
 11,~\cite{kit}). From~\cite[Theorem 3]{kit}, we know that
 \begin{equation}
   |A_n|=(n-1)!!+(n-2)!!,
 \end{equation}
 where $n!!$ is defined by $0!!=1$, and for $n\geq 1$
 \begin{equation}
n!!=\left\{\begin{array}{cc} n\cdot (n-2)\ldots 3\cdot 1 & \text{ if
} n \text{ is odd,} \\ n\cdot (n-2)\ldots 4\cdot 2 & \text{ if } n
\text{ is even.}
\end{array}\right.
 \end{equation}

 For $p\in\{213,231,312\}$ and $n\in\N$, let $\mathbf{p}_n$ denote the
 total number of occurrences of $\mathbf{p}$ in the permutations of
 $\Av_n(123,132,321)$. For $n \ge 2$, there are $n-2$ occurrences of a
 length 3 pattern in each permutation, so we have
 $\mathbf{213}_n+\mathbf{231}_n+\mathbf{312}_n=(n-2)|A_n|=(n-2)((n-1)!!+(n-2)!!).$ In the next theorem we give closed forms for $\mathbf{213}_n, \mathbf{231}_n,$ and $\mathbf{312}_n$.

 \begin{theorem}
   \begin{equation}
     \mathbf{231}_n=(n-1)!!\left\lceil\frac{n-3}{2}\right\rceil + (n-2)!!\left\lceil\frac{n-2}{2}\right\rceil,
   \end{equation}
   \begin{equation}
     \small\mathbf{312}_n=(n-1)!!\left(\frac{(-1)^{n-1}+n-3}{4} + \frac{1}{2}\sum_{\substack{k=1\\k\ne n\bmod{2}}}^{n-1} \frac{1}{k} \right)+(n-2)!!\left(\frac{(-1)^n+n-4}{4} + \frac{1}{2}\sum_{\substack{k=1\\k=n\bmod{2}}}^{n-2} \frac{1}{k} \right),
   \end{equation}
   \begin{equation}\mathbf{213}_n=(n-2)((n-1)!!+(n-2)!!)-\mathbf{231}_n-\mathbf{312}_n.
   \end{equation}
\end{theorem}

\begin{proof}
    From~\cite[Theorem 3]{kit}, we know that a permutation $\pi\in A_n$ is \textit{alternating} or \textit{reverse alternating}, namely $\pi=a_1\ldots a_n$ with $a_1>a_2<a_3\ldots$ or $a_1<a_2>a_3\ldots$. Moreover, for such a permutation we have either $a_n=1$, or $a_{n-1}=1$, and if we begin with either $a_n$ or $a_{n-1}$, whichever is equal to one, and build a sequence by including every other element as we read to the left, we get an increasing sequence. Let $A_n^r$ (resp. $A_n^l$) be the subset of $A_n$ consisting of permutations $\pi$ such that $a_n=1$ (resp. $a_{n-1}=1$). Again from the proof of~\cite[Theorem 3]{kit}, we have $|A_n^r|=(n-2)!!$ and $|A_n^l|=(n-1)!!$. From the structure of these permutations, it is easy to see that the positions of the occurrences of $231$ are exactly the indices $n-2,n-4,\ldots,2 \text{ or } 1$ for the permutations of $A_n^r$, and the indices $n-3,n-5,\ldots,2 \text{ or } 1$ for the permutations of $A_n^l$. The result on $\mathbf{231}_n$ follows.
    
    Now let us consider $\mathbf{312}_n$. Let $\mathbf{312}^r_n$ (resp. $\mathbf{312}^l_n$) denote the total number of occurrences of the consecutive pattern $312$ in $A_n^r$ (resp. $A_n^l$).
\begin{lemma}\label{rec 312}
    For any $n\geq 5$,
    \begin{enumerate}
        \item $\mathbf{312}^r_n=\mathbf{312}_{n-1}^l$,
        \item $\mathbf{312}^l_n=(n-1)(\mathbf{312}_{n-2}^l+(n-3)!!)-(n-5)!!\frac{(n-3)(n-2)}{2}.$
    \end{enumerate}
\end{lemma}
\begin{proof}
(1) follows from the fact that if $\pi\in A_n^r$, then the permutation induced by the $n-1$ first letters of $\pi$ is in $A_{n-1}^l$ (otherwise we would have an occurrence of $321$ on the last 3 letters of $\pi$). For (2) consider for $2\leq k\leq n$ the subset $B_n^k\subset A_n^l$ consisting of the permutations of $A_n^l$ with $a_n=k$. First assume that $k>2$, and let $\pi\in B_n^k$. Then $\pi=a_1\ldots a_{n-4}2a_{n-2}1k$ with $a_1\ldots a_{n-4},a_{n-2}\in\{3,\ldots,n\}\backslash\{k\}$. The permutation $\pi':=a_1\ldots a_{n-4}2a_{n-2}$ belongs to $A_{n-2}^l$. We know that $\pi\longrightarrow\pi'$ is a bijection from $B_n^k$ to $A_{n-2}^l$ that preserves the number of occurrences of $312$ on the first $n-2$ letters. It remains to count the number of occurrences of $312$ in the last 4 letters of each permutation of $B_n^k$, that is in $2a_{n-2}1k$ using our notation. Note that $2a_{n-2}1$ can never be an occurrence of $312$, and $a_{n-2}1k$ is one if and only if $a_{n-2}>k$. Therefore, each permutation of $A_{n-2}^l$ will have an additional occurrence of 312 in the corresponding permutation in $B_n^k$ unless the permutation ends in $2b1k$ with $3\leq b < k$. Since there are $(k-3)|A_{n-4}^l|$ such permutations, the total number of occurrences of $312$ in $B_n^k$ is $\mathbf{312}_{n-2}^l+|A_{n-2}^l|-(k-3)|A_{n-4}^l|$. Similarly, when $k=2$, permutations from $B_n^2$ have the form $a_1\ldots a_{n-2}12$, with $a_1\ldots a_{n-2}\in A_{n-2}^l$, otherwise the permutation would end in $4312$ which contains an occurrence of 321. Since $a_{n-2}12$ is always an occurrence of $312$, the total number of occurrences of $312$ in $B_n^2$ is then $\mathbf{312}_{n-2}^l+|A_{n-2}^l|$. We finally deduce that
\begin{align}
    \mathbf{312}_n^l &= \mathbf{312}_{n-2}^l+|A_{n-2}^l|+\sum_{k=3}^n (\mathbf{312}_{n-2}^l+|A_{n-2}^l|-(k-3)|A_{n-4}^l|)\\
    &=(n-1)(\mathbf{312}_{n-2}^l+(n-3)!!)-(n-5)!!\frac{(n-3)(n-2)}{2}.
\end{align}
\end{proof}

Let $u_n:=\frac{\mathbf{312}^l_n}{(n-1)!!}$ and $f(z):=\sum_{n=3}^\infty u_nz^n$. From Lemma~\ref{rec 312} (2) we deduce that $u_n=u_{n-2}+1-\frac{n-2}{2(n-1)}$ for any $n\geq 5$. We easily derive that
\begin{equation}f(z)=\frac{z(2(z-1)\ln(1-z)+z^3+3z^2-2z)}{4(1-z)^2(1+z)},
\end{equation}
and finally
\begin{equation}\mathbf{312}_n^l=(n-1)!!\left(\frac{(-1)^{n-1}+n-3}{4} + \frac{1}{2}\sum_{\substack{k=1\\k\ne n\bmod{2}}}^{n-1} \frac{1}{k} \right),
\end{equation}
and with Lemma~\ref{rec 312} (1),
\begin{equation}\mathbf{312}_n^r=(n-2)!!\left(\frac{(-1)^n+n-4}{4} + \frac{1}{2}\sum_{\substack{k=1\\k=n\bmod{2}}}^{n-2} \frac{1}{k} \right).
\end{equation}
\end{proof}

\begin{corollary}
    $\pop_{11}(231)=1/2$, and $\pop_{11}(213)=\pop_{11}(312)=1/4.$
\end{corollary}

\section{Avoiding $123$, $132$}
\label{sec:class17}

\cite{claes} proved that the Foata transform induces a bijection between $\Av_n(123,132)$ and $\mathcal{I}_n$, the set of involutions of length $n$, showing in particular that $|\Av_n(123,132)|=|\mathcal{I}_n|$. Let us recall briefly this process. An involution $\pi\in\mathcal{I}_n$ has all cycles of length 1 or 2. We introduce a standard form for writing $\pi$:
\begin{enumerate}
    \item Each cycle is written with its least element first.
    \item The cycles are written in decreasing order of their least element.
\end{enumerate}
Denote by $\Hat{\pi}$ the permutation obtained from $\pi$ by writing it in standard form and erasing the parentheses separating the cycles. Then $\pi\mapsto\Hat{\pi}$ is a bijection between $\mathcal{I}_n$ and $\Av_n(123,132)$. We will use this bijection to study the frequencies of each pattern of length 3 in $\Av_n(123,132)$, taking advantage of some knowledge we have on the involutions. Table~\ref{corr Av(123,132) and I_n} shows the correspondence between patterns in $\Av_n(123,132)$ and in $\mathcal{I}_n$.

\begin{table}[h]
  \begin{center}
    \begin{tabular}{|c|c|c|}
    \hline
        Pattern in $\Av_n(123,132)$ & Pattern in $\mathcal{I}_n$, with $a<b<c$ & Fixed point-free pattern in $\mathcal{I}_n$\\
        \hline
       $321$  &  $(c)(b)(a)$ or $(c)(b)(a \; \star)$ or $(\star\; c)(b)(a)$ & $\emptyset$\\
       \hline
        $231$ & $(b\; c)(a)$ or $(b\; c)(a\;\star)$ & $(b\; c)(a\;\star)$\\
        \hline
        $213$ & $(b)(a\; c)$ or $(\star\; b)(a\; c)$ & $(\star\; b)(a\; c)$\\
        \hline
        $312$ & $(c)(a\; b)$ or $(\star\; c)(a\; b)$ & $(\star\; c)(a\; b)$\\
        \hline
    \end{tabular}
    \caption{The correspondence between patterns in $\Av_n(123,132)$ and $\mathcal{I}_n$.}
    \label{corr Av(123,132) and I_n}
  \end{center}
\end{table}

\begin{example}
    The involution $\pi=732458169\in\mathcal{I}_9$ has standard form $\pi=(9)(6\; 8)(5)(4)(2\; 3)(1\; 7)$, so $\Hat{\pi}=968542317\in\Av_9(123,132)$. The occurrence 854 of the pattern $321$ in $\Hat{\pi}$ corresponds to the occurrence $(6\; 8)(5)(4)$ of the pattern $(\star\; c)(b)(a)$ in $\pi$.
\end{example}

Thanks to the Foata transform, we converted our problem of counting
patterns in $\Av_n(123,132)$ into counting patterns in
$\mathcal{I}_n$. In Lemma \ref{equ fixed point-free patterns in I_n}
we will see that we can actually focus on fixed point-free patterns in
$\mathcal{I}_n$. To do so, we first need the following lemma that
states that fixed points are rare within involutions.

\begin{lemma}\label{fixed points I_n}
  Let $\textsc{fp}_n$ be the total number of fixed points in $\mathcal{I}_n$. Then
  \begin{equation}
    \frac{\textsc{fp}_n}{|\mathcal{I}_n|}\underset{n\to\infty}{\sim}\sqrt{n}.
  \end{equation}
\end{lemma}
\begin{proof}
  The generating function $I(z,u)$ such that the coefficient of $z^nu^k$ is the number of involutions in $\mathcal{I}_n$ having $k$ fixed points, divided by $n!$, is given by $I(z,u)=e^{zu+z^2/2}.$ The average number of fixed points in $\mathcal{I}_n$ is then
  \begin{equation}
    \frac{\textsc{fp}_n}{|\mathcal{I}_n|}=\frac{[z^n]\partial_uI(z,u)\vert_{u=1}}{[z^n]I(z,1)}=\frac{[z^{n-1}]I(z,1)}{[z^n]I(z,1)}.
  \end{equation}

It is known, see for instance~\cite{flaj} and~\cite{wilf}, that 
\begin{equation}
  [z^n]I(z,1)\underset{n\to\infty}{\sim} \frac{n^{-n/2}}{2\sqrt{\pi
      n}}\exp\left(\frac{n}{2}+\sqrt{n}-\frac{1}{4}\right),
\end{equation}
which allows us to conclude after simplification\footnote{See also
Michael Lugo's blog post \\
\hspace*{2em}\url{http://godplaysdice.blogspot.com/2008/02/how-many-fixed-points-do-involutions.html}}.
\end{proof}

Lemma~\ref{fixed points I_n} indicates that fixed points are quite
rare within the involutions. Consequently, in order to compute the
number of occurrences of a pattern in $\Av_n(123,132)$, it suffices to
compute the number of occurrences of the corresponding fixed
point-free pattern in $\mathcal{I}_n$.

\begin{lemma}\label{equ fixed point-free patterns in I_n}
  Using the same notation as in Table \ref{corr Av(123,132) and I_n}, let $\boldsymbol{\alpha}:=(b\; c)(a\;\star)$, $\boldsymbol{\beta}:=(\star\; b)(a\; c)$ and $\boldsymbol{\gamma}:=(\star\; c)(a\; b)$ be the three consecutive fixed point-free patterns of length 3 in $\mathcal{I}_n$, and let $\boldsymbol{\alpha}_n$, $\boldsymbol{\beta}_n$ and $\boldsymbol{\gamma}_n$ denote their number of occurrences in $\mathcal{I}_n$, respectively. Then $\pop_{17}(321)=0$,
  \begin{equation}
    \pop_{17}(231)=\underset{n\to\infty}{\lim}\frac{\boldsymbol{\alpha}_n}{n|\mathcal{I}_n|},\quad \pop_{17}(213)=\underset{n\to\infty}{\lim}\frac{\boldsymbol{\beta}_n}{n|\mathcal{I}_n|},\quad \text{and}\quad \pop_{17}(312)=\underset{n\to\infty}{\lim}\frac{\boldsymbol{\gamma}_n}{n|\mathcal{I}_n|}.
  \end{equation}
\end{lemma}
\begin{proof}
    Let $\boldsymbol{\delta}_n$ be the total number of occurrences of a pattern admitting a fixed point in $\mathcal{I}_n$. In order to prove the lemma, it suffices to show that $\underset{n\to\infty}{\lim} \boldsymbol{\delta}_n/(n|\mathcal{I}_n|)=0$. Since each fixed point appears in 3 different occurrences of a pattern, we have $\boldsymbol{\delta}_n\leq 3 \cdot \textsc{fp}_n$. The result follows from Lemma~\ref{fixed points I_n}.
\end{proof}

\begin{proposition}\label{class 17 pop 321 213}
    $\pop_{17}(231)=1/2$.
\end{proposition}
\begin{proof}
    Lemma~\ref{equ fixed point-free patterns in I_n} indicates that it suffices to compute the popularity of the pattern $\boldsymbol{\alpha}$ among the fixed point-free patterns in $\mathcal{I}_n$. Such a pattern consists of 2 consecutive transpositions. But for each pair of consecutive transpositions in a permutation of $\mathcal{I}_n$, we have exactly one occurrence of the pattern $\boldsymbol{\alpha}$, and one occurrence of $\boldsymbol{\beta}$ or $\boldsymbol{\gamma}$. Thus, the popularity of $\boldsymbol{\alpha}$ among fixed point-free patterns in $\mathcal{I}_n$ is $\frac{1}{2}$, which finishes the proof.
\end{proof}

We now consider the pattern 213. By Lemma \ref{equ fixed point-free patterns in I_n}, it suffices to estimate the popularity of the pattern $(\star\; b)(a\; c)$ in the involutions. By definition of the standard form of the permutation, such an occurrence is necessarily of the form $(b\; c)(a\; d)$, with $a<b<c<d$. In $\Av(123,132)$, this corresponds to an occurrence of the pattern $2314$. Therefore, we compute the exponential generating function of $\mathbf{2314}_n$, the number of occurrences of 2314 in $\Av_n(123,132)$.

\begin{lemma}\label{formula 2314 class 17}
  $\mathbf{2314}_4=1$, $\mathbf{2314}_5=4$, and for all $n\geq6$,
  \begin{equation}
    \mathbf{2314}_n=\mathbf{2314}_{n-1}+(n-1)\mathbf{2314}_{n-2}+\binom{n-2}{2}|\mathcal{I}_{n-4}|.
  \end{equation}
\end{lemma}
\begin{proof}
    It is easy to see that each permutation in Class 17 has 1 in
    either the first or the second position from the
    right. Let $A_n^r$ (resp. $A_n^l$) be the subset of $\Av_n(123,132)$ consisting of permutations $\pi=a_1\ldots a_n$ such that $a_n=1$ (resp. $a_{n-1}=1$). The number of occurrences of $2314$ in each permutation of $A_n^r$ is exactly $\mathbf{2314}_{n-1}$. Let
    $\pi\in A_n^l$, then $\pi=a_1\ldots a_{n-2}1k$
    for some $k\in\{2,\ldots,n\}$, and $a_1\ldots
    a_{n-2}\in\Av_{n-2}(123,132)$. So for each occurrence of
    $2314$ in a permutation of $\Av_{n-2}(123,132)$, we have
    $n-1$ occurrences of $2314$ in a permutation of
    $A_n^l$. However, occurrences can also appear in
    the last four positions of permutations from $A_n^l$: there is one occurrence for each permutation ending with
    $2i1k$, for $3\leq i < k\leq n$, and there are
    $\binom{n-2}{2}|\Av_{n-4}(123,132)|=\binom{n-2}{2}|\mathcal{I}_{n-4}|$ such permutations. In the end,
    we obtain the desired formula.
\end{proof}

\begin{proposition}
  Let $G(z)=\sum_{n=4}^{+\infty}\frac{\mathbf{2314}_n}{n!}z^n$ be the EGF of $(\mathbf{2314}_n)_{n\geq4}$. Then
  \begin{equation}
    G(z)=\frac{e^{\frac{(1+z)^2}{2}}}{2}\int_{0}^z e^{-\frac{(1+t)^2}{2}}dt+\frac{z(z-2)e^{z+\frac{z^2}{2}}}{4}.
  \end{equation}
\end{proposition}

\begin{proof}
  From Lemma~\ref{formula 2314 class 17} we can verify that $G$ satisfies the following Cauchy problem:
  \begin{equation}
    \left\{
    \begin{array}{l}
        G''(z)-(1+z)G'(z)-G(z)=\frac{z^2}{2}e^{z+\frac{z^2}{2}},  \\
        G(0)=G'(0)=0.
    \end{array}\right.
  \end{equation}
  Note that $e^{z+\frac{z^2}{2}}$ is the EGF of
  $(|\mathcal{I}_n|)_{n\geq0}$. Solving this differential equation, we
  get that its only solution is the one stated in the lemma.
\end{proof}

\begin{corollary}
    $\pop_{17}(312)=\pop_{17}(213)=1/4$.
\end{corollary}
\begin{proof}
  We proceed in two steps. First we prove that
  $[z^n]z(z-2)e^{z+\frac{z^2}{2}}\sim\frac{n|\mathcal{I}_n|}{n!}$, and
  secondly we show that $[z^n]\left(e^{\frac{(1+z)^2}{2}}\int_{0}^z
  e^{-\frac{(1+t)^2}{2}}dt\right)=o\left(\frac{n|\mathcal{I}_n|}{n!}\right)$. Those
  two facts prove that $[z^n]G(z)\sim\frac{n|\mathcal{I}_n|}{4\cdot
    n!}$, and so $\pop_{17}(2314)=\pop_{17}(213)=1/4$. By Proposition
  \ref{class 17 pop 321 213} we then have $\pop_{17}(312)=1/4$. For
  the first point, we already know (see for instance~\cite{flaj}
  and~\cite{wilf}) that
  \begin{equation}[z^n]e^{z+\frac{z^2}{2}}=\frac{|\mathcal{I}_n|}{n!}\underset{n\to\infty}{\sim} \frac{n^{-n/2}}{2\sqrt{\pi
        n}}\exp\left(\frac{n}{2}+\sqrt{n}-\frac{1}{4}\right).
  \end{equation}
  We deduce after simplification that
  \begin{equation}
    \frac{|\mathcal{I}_{n-1}|}{(n-1)!}\sim\frac{n^{-n/2}}{2\sqrt{\pi}}\exp\left(\frac{n}{2}+\sqrt{n}-\frac{1}{4}\right), \mbox{\ and \ } \frac{|\mathcal{I}_{n-2}|}{(n-2)!}\sim\frac{\sqrt{n}\cdot n^{-n/2}}{2\sqrt{\pi}}\exp\left(\frac{n}{2}+\sqrt{n}-\frac{1}{4}\right).
  \end{equation}
  Consequently,
  \begin{equation}[z^n]z(z-2)e^{z+\frac{z^2}{2}}=\frac{|\mathcal{I}_{n-2}|}{(n-2)!}-2\frac{|\mathcal{I}_{n-1}|}{(n-1)!}\sim\frac{\sqrt{n}\cdot n^{-n/2}}{2\sqrt{\pi}}\exp\left(\frac{n}{2}+\sqrt{n}-\frac{1}{4}\right)\sim\frac{n|\mathcal{I}_n|}{n!},
  \end{equation}
which proves the first point. For the second one, let us first consider the coefficients of the other term of $G(z)$.
\begin{lemma}\label{positivity coeffs F(z)}
    Let $F(z)=e^{\frac{(1+z)^2}{2}}\int_{0}^z e^{-\frac{(1+t)^2}{2}}dt$. Then for all $n\geq0$, $[z^n]F(z)\in\frac{1}{n!}\cdot\N$.
\end{lemma}

\begin{proof}
    $F$ satisfies the following equality: $F'(z)=(1+z)F(z)+1$. Then, with a direct induction, there exist two polynomials $p_n,q_n\in\N[z]$, with $\deg(p_n)=n$ and $\deg(q_n)=n-1$ such that $F^{(n)}(z)=p_n(z)F(z)+q_n(z)$. Since $F(0)=0$, we have $F^{(n)}(0)\in\N$, and as $[z^n]F(z)=\frac{F^{(n)}(0)}{n!}$, the result follows\footnote{$F'(z)$ is the EGF of the sequence
    \href{https://oeis.org/A000932}{A000932} in OEIS,~\cite{OEIS}.}.
\end{proof}

Lemma \ref{positivity coeffs F(z)} ensures the positivity of the coefficients of $F(z)$. We can then apply a saddle point bound to $F$ (see~\cite[Corollary VIII.1]{flaj}). Indeed, by the positivity, $F(z)z^{-n-1}$ has a unique saddle point $\zeta$ defined by \begin{equation}\label{eq: saddle point equation}
    \zeta\frac{F'(\zeta)}{F(\zeta)}=n+1, \mbox{ or equivalently } \frac{\zeta}{F(\zeta)}=n+1-\zeta-\zeta^2,
\end{equation} and then 
\begin{equation}\label{eq: saddle point bound}
    [z^n]F(z)\leq\frac{F(\zeta)}{\zeta^n}.
\end{equation}
Given the saddle point equation (\ref{eq: saddle point equation}), it seems out of reach to obtain an exact expression of $\zeta$. However, the saddle point is the value minimizing the right-hand term in (\ref{eq: saddle point bound}), and since we are only looking for a good enough upper bound on $[z^n]F(z)$, i.e. a bound which is $o\left(\frac{n|\mathcal{I}_n|}{n!}\right)$, a nice approximation of the saddle point may yield a sufficient bound. It turns out that choosing $\zeta=\sqrt{n}$ is sufficient. The upper bound (\ref{eq: saddle point bound}) then becomes
\begin{align}
    [z^n]F(z)&\leq \frac{e^{\frac{(1+\sqrt{n})^2}{2}}\int_0^{\sqrt{n}}e^{-\frac{(1+t)^2}{2}}dt}{(\sqrt{n})^n}\\
    &\leq \left(\int_0^{+\infty}e^{-t-\frac{t^2}{2}}dt\right)\cdot n^{-n/2}\cdot \exp\left(\frac{n}{2}+\sqrt{n}\right).
\end{align}

This is indeed enough to prove that $[z^n]F(z)=o\left(\sqrt{n}\cdot n^{-n/2}\cdot\exp\left(\frac{n}{2}+\sqrt{n}\right)\right)=o\left(\frac{n|\mathcal{I}_n|}{n!}\right)$.
\end{proof}

\section{Open questions}\label{sec:open}

The conjecture, presented below, is true for cases solved in our work.
Is it true in the general case?
\begin{conjecture}
  For $m>0$, denote by $p_1,\ldots,p_k,p$ certain consecutive patterns
  of length $m$. Let $\mathcal{A}_n=\Av_n(p_1,\ldots,p_k)$ and
  $\mathcal{B}_n=\Av_n(p_1,\ldots,p_k,p)$. Suppose that
  $\pop_\mathcal{A}(p)=0$, then for every consecutive pattern $q$ of
  length $m$ we have $\pop_\mathcal{B}(q)=\pop_\mathcal{A}(q)$.
\end{conjecture}

It may also be interesting to answer the following questions:
\begin{enumerate}
    \item Can we solve Classes 10, 12, 13, 14 and 15 using an
      enumerative or probabilistic approach?
      Our numerical experiments suggest that the convergence rate
      is too low to formulate conjectures
      in these cases.
    \item What happens when we avoid only one consecutive pattern of
      length 3?
    \item Can we find a set of patterns, such that when we avoid them we obtain
      irrational asymptotic popularity for some remaining pattern?
    \item Does the limit~\eqref{eq:limpop} always exist? 
    If not, can we characterize the patterns for which this limit exists?
\end{enumerate}

\acknowledgements

We would like to thank Sergi Elizalde for helpful conversations about
the patterns of this paper, and anonymous reviewers for their careful
reading, comments and advice. This research was funded, in part, by
the Agence Nationale de la Recherche (ANR), grants ANR-22-CE48-0002,
ANR-25-CE48-0602 and by the Regional Council of
Bourgogne-Franche-Comté.

\nocite{*}
\bibliographystyle{abbrvnat}
% use the following instead if you encounter problems 
%\bibliographystyle{alpha}
\bibliography{b}
\label{sec:biblio}

\end{document}